\documentclass{article}
\usepackage[utf8]{inputenc}
\usepackage{amssymb}
\usepackage{amsmath}
\title{The Weingarten map and curvatures of $SL(n, \mathbb{R})$} 
\author{Supravat Sarkar}
\date{}
\begin{document}

\maketitle
\begin{abstract}
    In this article, we shall derive by elementary calculations the Gauss map, spherical image, Weingarten map and the curvatures at identity of the special linear group, that is, the matrices of determinant 1. We could not find any reference where this simple, nice computation has been explicitly put down.
\end{abstract}
\begin{center}
\end{center}
\begin{center}
\end{center}

\section{Introduction}
\vskip 5mm

\noindent For a surface in $\mathbb{R}^3$, the Weingarten map at a
point describes the variation of the normal field with respect to a
tangent vector at the point. In other words, it gives an extrinsic
information on the shape of the surface. Formally, the Weingarten
map of a hypersurface in $\mathbb{R}^N$ can be defined as the
derivative of the Gauss map. Knowing the Weingarten map enables us
to compute the principal curvatures. Gaussian curvature and the mean curvature at
various points. In this note, we consider the group $SL(n,
\mathbb{R})$ as a hypersurface in $\mathbb{R}^{n^2}$ by means of the
determinant map and compute the Weingarten map and deduce the
curvatures etc. We could not find any reference where this simple,
nice computation has been explicitly put down. \vskip 5mm

\section{Notations}

\noindent
 Let $S$=$SL(n,$ $\mathbb{R})$=$f^{-1}(1)$, where
 $f:M_{n}(\mathbb{R})\to \mathbb{R}$ is given by $f(A)=det(A).$
It is easy to see that $\nabla{f}$ is nowhere $0$ in $S$, so $S$ is
an $n^{2}-1$-surface in $\mathbb{R}^{n^{2}}$, the latter identified
with $M_{n}(\mathbb{R})$. We orient $S$ by the unit normal vector
field $N=$ $\nabla{f}/||\nabla{f}||$. In this article, we compute
the Gauss map, spherical image of $S$, Weingarten map and finally the curvatures of
$S$ at $I$.\\

\noindent
   We mostly follow the notations of Thorpe's book (\cite{T}). In addition, we adopt the following convention: when we write a tangent vector to $S$ at a point $p\in S$ , we just write the vector part, do not mention the point $p$. In this convention, we give the the definitions of the relevant terms. All of them can be found in the corresponding chapters of (\cite{T}).
   \begin{enumerate}
   \item The tangent space $S_{p}$ of $S$ at a point $p\in S$ is the set of all tangent vectors to $S$ at $p$. So, this is a vector subspace of $R^{n^2}$.
       \item \textit{Gauss map} on $S$ is the map $N:S\to S^{n^2-1}$.
       \item \textit{Spherical image} of $S$ is the range of the Gauss map.
       \item \textit{Weingarten map} of $S$ at a point $p\in S$ is the map $L_{p}:S_{p}\to S_{p}$ defined by $L_{p}(v)=-(N\circ \alpha)^{\prime}(0)$, where $\alpha:I\to S$ is a smooth curve with $\alpha(0)=p$, $\alpha^{\prime}(0)=v$. (I is an open interval in $\mathbb{R}$ containing $0$.) It can be checked that $L_{p}(v)$ is independent of the choice of $\alpha$ and $L_{p}(v)=-\widetilde{N}^{\prime}(p)(v)$, where
       $\widetilde{N}:U\to R^{n^2}$ is an extension of the Gauss map $N$ to an open set $U$ of $R^{n^2}$ containing $p$.
       \item \textit{Tangent space} at $I$ of $S$ can be identified with
       the space $sl_n$ of trace $0$ matrices. Further, $sl_{n}=(sym_{n}\cap sl_{n})\oplus sk_{n},$
where $sym_{n}, sk_{n}$ are the subspaces of symmetric and
skew-symmetric matrices of $M_{n}(\mathbb{R})$, respectively.
   \end{enumerate}

\vskip 5mm

   \section{Results}
   We shall prove:
   \vskip 2mm

\noindent {\bf Theorem.} For the group $S=SL(n, \mathbb{R})=
f^{-1}(1)$, where $f:M_{n}(\mathbb{R})\to \mathbb{R}$ is given by
$f(A)=det(A),$ we have:
   \begin{enumerate}
       \item The Gauss map of $S$ is given by $N(A)=(A^{-1})^{T}/||(A^{-1})||.$
       \item The spherical image of $S$ is $\{A \in S^{n^2-1}: det(A)>0\}$. So, the spherical image is symmetric about $0$ if and only if $n$ is even.
       \item The Weingarten map of $S$ at $I$ is given by
       $L_{I}(H)=n^{-1/2} H^{T}$, $\forall H\in S_{I}$.
       \item The principal curvatures of $S$ at $I$ are:\\
   $n^{-1/2}$ and $-n^{-1/2}$ with respective multiplicities
   $\frac{n^{2}+n-2}{2}$ and $\frac{n^{2}-n}{2}.$
   \end{enumerate}

\vskip 5mm

 \section{Gauss map and spherical image}
  The Gauss map $N$ of $f$ would be a $n\times n$ matrix.  $\nabla{f}(A)_{ij}$=$\frac{\delta f}{\delta a_{ij}} (A)$=$cof_{A}(i,j)$, as by Laplace expansion of determinant, $f(A)=\sum_{j=1}^{n} a_{ij}cof_{A}(i,j)$, and none of these cofactors involve any $a_{ij}$. The $n\times n$ matrix $B$ with $(i,j)$'th entry $\nabla{f}(A)_{ij}$ is then the transpose of the adjugate matrix of $A$, but since determinant of $A$ is 1, its adjugate is same as its inverse. So $B$=$(A^{-1})^{T}$. The Gauss map $N:S\to S^{n^2-1}$ is given by $N(A)=B/||B||=(A^{-1})^{T}/||A^{-1}||.$ Here for a matrix $C=((c_{ij})),$ $||C||=\sqrt{\sum_{i,j} c_{ij}^{2}}$.

  \noindent The spherical image of $S$ is given by
$$\{(A^{-1})^{T}/||A^{-1}||: det(A)=1\}$$
$$=\{A/||A||: det(A)=1\}$$
$$=\{A \in S^{n^2-1}: det(A)>0\},$$ because det($A/||A||$)=$||A||^{-n}>0$,
and because $det(d^{-1/n} A)=1$, and $A=(d^{-1/n} A)/||d^{-1/n} A||$
if $A \in S^{n^2-1}$ is such that $det(A)=d>0$.

\vskip 5mm

  \section{Weingarten map}
  Tangent space of S at I is $S_{I}$=$\{H\in M_{n}(\mathbb{R}):$ and $tr(H)=0\}=sl_{n},$ as it is the orthogonal complement of $N(I)$.
  The Weingarten map is
  $$L_{I}(H)=-\widetilde{N}^{\prime}(I)(H),$$
  where $\widetilde{N}$ is an extension of the Gauss map $N$. We may take $\widetilde{N}(A)=(A^{-1})^{T}/||A^{-1}||$ on the open set $U$ of
  invertible matrices.\\
  For $1\leq i,j \leq n$, let $E_{ij}$ denote the $n\times n$ matrix whose $(i,j)$'th entry is 1 and all other entries are
  $0$.\\
  For $i\neq j,$
  $$\widetilde{N}^{\prime}(I)(E_{ij})=\frac{d}{dt} \Bigr\rvert_{t = 0} \widetilde{N}(I+tE_{ij})=\frac{d}{dt}
  \Bigr\rvert_{t = 0}(\frac{I-tE_{ji}}{\sqrt{n+t^{2}}})=-n^{-1/2}E_{ji},$$
  as $$\frac{d}{dt} \Bigr\rvert_{t =
0}\frac{1}{\sqrt{n+t^{2}}}=-\frac{1}{2.(n+t^{2})^{3/2}}.2t\Bigr\rvert_{t
= 0}=0$$ and $$\frac{d}{dt} \Bigr\rvert_{t =
0}\frac{-t}{\sqrt{n+t^{2}}}=-n^{-1/2}.$$

\noindent Hence
  $\widetilde{N}^{\prime}(I)(E_{ii})=\frac{d}{dt} \Bigr\rvert_{t = 0}
  \widetilde{N}(I+tE_{ij})$.\\
  Now,
$$\widetilde{N}(I+tE_{ii})= \frac{1}{\sqrt{n-1+(1+t)^{-2}}}
(I+((1+t)^{-1}-1)E_{ii})$$ $$=\frac{1}{\sqrt{(n-1)(1+t)^{2}+1}}
  ((1+t)I-tE_{ii}).$$

  \noindent
  Since
  $$\frac{d}{dt} \Bigr\rvert_{t = 0} \frac{1+t}{\sqrt{(n-1)(1+t)^{2}+1}}=\frac{\sqrt{n}-\frac{1}{2\sqrt{n}}
  2(n-1)}{n-1+1}=\frac{2}{2n\sqrt{n}}=n^{-3/2}$$
  and
$$\frac{d}{dt} \Bigr\rvert_{t = 0}
\frac{1}{\sqrt{(n-1)(1+t)^{2}+1}}=-\frac{1}{2}
  \frac{2(n-1)}{n^{3/2}}=-\frac{n-1}{n^{3/2}},$$
  we obtain
  $$\widetilde{N}^{\prime}(I)(E_{ii})=n^{-3/2}I-n^{-1/2}E_{ii}.$$
   Thus, what we have got is that
   $$\widetilde{N}^{\prime}(I)(E_{ij})= -n^{-1/2}E_{ji} ~\forall~ i\neq j$$ and
   $$\widetilde{N}^{\prime}(I)(E_{ii})= n^{-3/2}I-n^{-1/2}E_{ii}.$$
From the above values, we obtain
   $$\widetilde{N}^{\prime}(I)(E_{ii}-E_{jj})=-n^{-1/2}(E_{ii}-E_{jj}).$$
   As any trace zero matrix is a real linear combination of the matrices $E_{ii}-E_{jj}$ and $E_{ij}(i\neq j)$ ,
   we see that $\widetilde{N}^{\prime}(I)(H)=n^{-1/2}H^{T}$ $\forall H\in S_{I}.$ \\
   Hence, the Weingarten map is given by
   $L_{I}(H)=n^{-1/2} H^{T}$, $\forall H\in S_{I}$.

\vskip 5mm

   \section{Calculation of curvatures}

  \noindent We need to find the Eigenvalues of $L_{I}$. We know that
$sl_{n}=(sym_{n}\cap sl_{n})\oplus sk_{n}$, where $sym_{n}, sk_{n}$
are the subspaces of symmetric and skew-symmetric matrices of
$M_{n}(\mathbb{R})$, respectively.\\

\noindent  So,
$$L_{I}(H)=n^{-1/2} H~~\forall~~l H\in sym_{n}\cap
sl_{n},$$ and
$$L_{I}(H)=-n^{-1/2} H~~\forall~~ H\in sk_{n}$$
where the dimension of $sym_{n}\cap sl_{n}$ is
$\frac{n^{2}+n-2}{2}$, and the dimension of $sk_{n}$ is
$\frac{n^{2}-n}{2}.$\\
Hence, we obtain that the principal curvatures of $S$ at
$I$ are:\\
$n^{-1/2}$ with multiplicity $\frac{n^{2}+n-2}{2}$, and $-n^{-1/2}$
with multiplicity $\frac{n^{2}-n}{2}.$\\
The Gauss-Kronecker curvature at $I$ is $$(-1)^{\frac{n^{2}-n}{2}}
   n^{-\frac{n^{2}-1}{2}},$$ and the mean curvature at $I$ is $$\frac{1}{n^{2}-1} \frac{1}{\sqrt{n}}
(\frac{n^{2}+n-2}{2}-\frac{n^2-n}{2})=\frac{1}{\sqrt{n}(n+1)}.$$
Finally, the 1st fundamental form at $I$ is $H\to tr(H^{T}H)$, 2nd
fundamental form is $H\to \frac{1}{\sqrt{n}} tr(H^{2}),$ $\forall$
$H\in sl_{n}$, as seen easily.
\section{Acknowledgement}
 I am thankful to Prof. Aniruddha Naolekar for encouragement and motivation.

\end{document}